\documentclass{article}

\usepackage{
 amsmath,
 amsxtra,
 amsthm,
 amssymb,
 mathrsfs
  }
\usepackage[all]{xy}

\newtheorem{theorem}{Theorem}[section]
\newtheorem{lemma}[theorem]{Lemma}

\newtheorem{proposition}[theorem]{Proposition}
\newtheorem{corollary}[theorem]{Corollary}
\newtheorem*{thm2}{Theorem}
\newtheorem*{cor2}{Corollary}

\theoremstyle{definition}
\newtheorem{defn}[theorem]{Definition}
\newtheorem{remark}[theorem]{Remark}

\newcommand{\bd}{\begin{defn}}
\newcommand{\ed}{\end{defn}}
\newcommand{\bl}{\begin{lemma}}
\newcommand{\el}{\end{lemma}}
\newcommand{\bp}{\begin{proposition}}
\newcommand{\ep}{\end{proposition}}
\newcommand{\bt}{\begin{theorem}}
\newcommand{\et}{\end{theorem}}
\newcommand{\bc}{\begin{corollary}}
\newcommand{\ec}{\end{corollary}}
\newcommand{\br}{\begin{remark}}
\newcommand{\er}{\end{remark}}
\newcommand{\ba}{\begin{array}}
\newcommand{\ea}{\end{array}}
\newcommand{\bpf}{\begin{proof}}
\newcommand{\epf}{\end{proof}}

\newcommand{\Q}{\mathbb{Q}}
\newcommand{\Qp}{\mathbb{Q}_p}

\newcommand{\Zp}{\mathbb{Z}_p}
\newcommand{\Op}{\mathcal{O}}
\newcommand{\Ga}{\Gamma}
\newcommand{\ga}{\gamma}
\newcommand{\La}{\Lambda}

\DeclareMathOperator{\Gal}{Gal}
\DeclareMathOperator{\Hom}{Hom}

\DeclareMathOperator{\Ext}{Ext}
\newcommand{\Iw}{\mathrm{Iw}}

\newcommand{\Tr}{\mathrm{Tr}}

\newcommand{\lra}{\longrightarrow}

\newcommand{\ot}{\otimes}

\newcommand{\ps}[1]{[[ #1 ]]}

\newcommand{\ilim}{\displaystyle \mathop{\varinjlim}\limits}
\newcommand{\plim}{\displaystyle \mathop{\varprojlim}\limits}

\numberwithin{equation}{section}

\begin{document}
\title{Comparing direct limit and inverse limit of even $K$-groups in multiple $\Zp$-extensions}
 \author{
  Meng Fai Lim\footnote{School of Mathematics and Statistics $\&$ Hubei Key Laboratory of Mathematical Sciences,
Central China Normal University, Wuhan, 430079, P.R.China.
 E-mail: \texttt{limmf@ccnu.edu.cn}} }
\date{}
\maketitle

\begin{abstract} \footnotesize
\noindent Iwasawa first established a duality relating the direct limit and the inverse limit of class groups in a $\Zp$-extension, and this result has recently been extended to multiple $\Zp$-extensions by many authors. In this paper, we establish an analogous duality for the direct limit and the inverse limit of higher even $K$-groups in a $\Zp^d$-extension.


\end{abstract}

\section{Introduction}

Let $F$ be a number field and let $K$ be a $\Zp$-extension of $F$. For every finite extension $L$ of $F$ contained in $K$, we let $A_L$ denote the $p$-primary part of the ideal class group of $L$. The ring of integer of $L$ is then denoted by $\Op_L$. Suppose that $L'$ is an extension of $L$ contained in $K$. Then there is a natural map $A_L \lra A_{L'}$ induced by the natural inclusion $\Op_L\subseteq \Op_{L'}$. On the other hand, we  have a map $A_{L'} \lra A_{L}$
going the other way which is induced by the (ideal) norm. Then one has
\[  \ilim_L A_L \quad \mbox{and}\quad \plim_L A_L,\]
where the direct limit (resp., the inverse limit) is taken with respect to the inclusion maps (resp., the norm maps). These two limit modules come naturally equipped with  $\Zp\ps{\Ga}$-module structures, where $\Ga= \Gal(K/F)\cong\Zp$.
For a $\Zp\ps{\Ga}$-module $M$, we write $M^\iota$ for the $\Zp\ps{\Ga}$-module which is the same underlying $\Zp$-module $M$ but whose $\Ga$-action is given by
\[ \ga \cdot_{\iota} x = \ga^{-1}x, \quad \gamma \in\Ga, x\in M.\]
A classical theorem of Iwasawa \cite{Iw73} then asserts that there is a pseudo-isomorphism
\[  \left(\plim_L A_L\right)^\iota \sim \left(\ilim_L A_L\right)^\vee \]
of $\Zp\ps{\Ga}$-modules. Here $(-)^\vee$ is the Pontryagin dual. This result of Iwasawa has been generalized to the context of a $\Zp^d$-extension (see the works of Nekov\'a\v{r} \cite{Ne}, Vauclair \cite{Vau} and, more recently, that of Lai and Tan \cite{LT}).

In this paper, we will consider the situation of the higher even $K$-groups. As before, let $p$ be a prime, and let $F$ be a number field. In the event that $p=2$, we shall assume further that the number field $F$ has no real primes. Let $F_\infty$ be a $\Zp^d$-extension of $F$. Fix an integer $i\geq 2$. For each finite intermediate extension $L$ of $F_\infty/F$, the works of Quillen \cite{Qui73b} and Borel \cite{Bo} tell us that the higher even $K$-group $K_{2i-2}(\Op_L)$ is finite. As is standard in Iwasawa theory, we are interested in the Sylow $p$-subgroup $K_{2i-2}(\Op_L)[p^\infty]$ of $K_{2i-2}(\Op_L)$. Now, for two finite subextensions $L\subseteq L'$, the inclusion $\Op_L\lra \Op_L'$ induces a map $\jmath_{L/L'}: K_{2i-2}(\Op_L)[p^\infty]\lra K_{2i-2}(\Op_{L'})[p^\infty]$ by functoriality. In the other direction, there is the norm map (also called the transfer map) $\Tr_{L'/L}: K_{2i-2}(\Op_{L'})[p^\infty]\lra K_{2i-2}(\Op_{L})[p^\infty]$. Similar to the class groups situation, we consider the following direct limit and inverse limit
\[  \ilim_L K_{2i-2}(\Op_L)[p^\infty] \quad \mbox{and}\quad \plim_L K_{2i-2}(\Op_L)[p^\infty],\]
whose transition maps are given by the maps $\jmath_{L/L'}$ and $\Tr_{L'/L}$ respectively. Again, these limit modules come equipped with natural $\Zp\ps{G}$-module structures, where $G=\Gal(F_\infty/F)\cong \Zp^d$. For a $\Zp\ps{G}$-module $M$, the module $M^\iota$ is defined similarly as before. The main result of this paper is then as follows.

\begin{thm2}[Theorem \ref{main theorem}]
Retain the notation as above. Then there is a pseudo-isomorphism
\[  \left(\plim_L K_{2i-2}(\Op_L)[p^\infty]\right)^\iota \sim\left(\ilim_L K_{2i-2}(\Op_L)[p^\infty]\right)^\vee\]
of $\Zp\ps{G}$-modules.
\end{thm2}

Since $K_0(\Op_L)[p^\infty] = A_L$, our result may therefore be interpreted as a generalization of the previous results of Iwasawa \textit{et al} to the higher even $K$-groups. For the remainder of the introductional section, we shall sketch the ideas of our proof, leaving the details to the body of the paper.

Let $M$ be a finitely generated $\Zp\ps{G}$-module. For two open subgroups $V\subseteq U$ of $G$, the norm map $N_{U/V}$ on $M_V$ factors through $M_U=(M_V)_{U/V}$ to yield a map $M_U\lra M_V$, which by abuse of notation is also denoted by $N_{U/V}$. Then $\{M_U\}$ forms a direct system.
We then say that the $\Zp\ps{G}$-module $M$ is \textit{systematically coinvariant-finite} if $M_U$ is finite for every open subgroup $U$ of $G$. For such a systematically coinvariant-finite module $M$, the direct limit $\ilim_U M_U$ is naturally a discrete $\Zp\ps{G}$-module. The following algebraic result, which will be proved in Section \ref{Iwasawa modules}, is an important ingredient for the proof of our main theorem.

\begin{thm2}[Theorem \ref{alg main}]
Let $M$ be a finitely generated $\Zp\ps{G}$-module which is systematically coinvariant-finite. Then we have an isomorphism
\begin{equation} \label{alg iso} \Ext^1_{\Zp\ps{G}}(M,\Zp\ps{G})\cong \Big( \ilim_U M_U\Big)^\vee \end{equation}
of $\Zp\ps{G}$-modules.
\end{thm2}


Since a systematically coinvariant-finite $\Zp\ps{G}$-module is automatically torsion over $\Zp\ps{G}$ (see Lemma \ref{finite fg tor} below), the left hand of the isomorphism (\ref{alg iso}) is pseudo-isomorphic to $M^{\iota}$ (for instance, see \cite[Proposition 8]{PR}). In view of this, to prove our Theorem \ref{main theorem}, it suffices to take $M =  \plim_L K_{2i-2}(\Op_L)[p^\infty]$, and show that this module satisfies the following properties:
\begin{enumerate}
  \item[(I)] For every open subgroup $U$ of $G$, we have an isomorphism
      \[t_U: \left(\plim_L K_{2i-2}(\Op_L)[p^\infty]\right)_U \cong K_{2i-2}(\Op_{L_U})[p^\infty],\] where $L_U$ is the fixed field of $U$. In particular, the $\Zp\ps{G}$-module $\plim_L K_{2i-2}(\Op_L)[p^\infty]$ is systematically coinvariant-finite.
  \item[(II)] For open subgroups $V\subseteq U$ of $G$, the following diagram \begin{equation} \label{K commute intro} \entrymodifiers={!! <0pt, .8ex>+} \SelectTips{eu}{}\xymatrix{
     \left(\plim_L K_{2i-2}(\Op_L)[p^\infty]\right)_U \ar[r]_(0.55){\sim}\ar[d]^{N_{U/V}} &  K_{2i-2}(\Op_{L_U})[p^\infty] \ar[d]^{\jmath_{L_U/L_V}}\\
     \left(\plim_L K_{2i-2}(\Op_L)[p^\infty]\right)_V \ar[r]_(0.55){\sim} &  K_{2i-2}(\Op_{L_V})[p^\infty]
     } \end{equation}
     commutes.
\end{enumerate}
The verification of the above two properties will be dealt with in Section \ref{Arithmetic preliminaries} (in particular, see Proposition \ref{main prop}).
As a corollary, we have the following result which is reminiscent to that for the class groups  (see \cite[Proposition 2]{Ba}, \cite[Corollary 1.2.1]{LT} and \cite[Th\'{e}or\`{e}me 4.1 and Th\'{e}or\`{e}me 4.4]{LaNQD}).

\begin{cor2}[Corollary \ref{pseudo K=0}]
Retain the above notation. Then $\plim_L K_{2i-2}(\Op_L)[p^\infty]$ is pseudo-null over $\Zp\ps{G}$ if and only if $\ilim_L K_{2i-2}(\Op_L)[p^\infty]=0$.
\end{cor2}

\subsection*{Acknowledgement}
This research is supported by the
National Natural Science Foundation of China under Grant No. 11771164.

\section{Systematically coinvariant-finite modules} \label{Iwasawa modules}

Throughout this section, $p$ will always denote a fixed prime. Let $d\geq 1$. We shall denote by $G$ the group isomorphic to the $d$-copies of the additive group $\Zp$. The completed group algebra $\Zp\ps{G}$ is defined by
\[ \plim_{U}\Zp[G/U],\]
where $U$ runs through all open subgroups of $G$ and the transition maps are given by the natural projection $\Zp[G/V]\twoheadrightarrow \Zp[G/U]$ for $V\subseteq U$. It is well-known that the ring $\Zp\ps{G}$ can be identified with the power series ring in $d$ variables over $\Zp$. In particular, it is a local ring.

\bd
Let $M$ be a $\Zp\ps{G}$-module. For an open subgroup $U$ of $G$, we write $M_U$ for the largest quotient
of $M$ on which $U$ acts trivially. The $\Zp\ps{G}$-module $M$ is then said to be \textit{systematically coinvariant-finite} if $M_U$ is finite for every open subgroup $U$ of $G$. As seen in the introduction, $\{M_U\}_U$ forms a direct system of finite modules with transition maps given by \[ N_{U/V}: M_U\lra M_V\]
for $V\subseteq U$. In particular, $\ilim_U M_U$ is a discrete $\Zp\ps{G}$-module.
\ed

\bl \label{finite fg tor}
A systematically coinvariant-finite $\Zp\ps{G}$-module $M$ is finitely generated torsion over $\Zp\ps{G}$.
\el

\bpf
Let $\mathfrak{m}$ be the (unique) maximal ideal of $\Zp\ps{G}$. Then it contains the augmentation ideal $I_G$ (for instances, see \cite[Proposition 5.2.16]{NSW}). Thus, $M/\mathfrak{m}M$ is a quotient of $M_G$ and hence finite. By the Nakayama lemma, this in turn implies that $M$ is finitely generated over $\Zp\ps{G}$. Finally, since $G\cong\Zp^d$, it is in particular a solvable $p$-adic Lie group. Therefore, we may apply the main result of \cite{BH} to conclude that $M$ is torsion over $\Zp\ps{G}$.
\epf

We can now state the main theorem of this section.

\bt \label{alg main}
Let $M$ be a systematically coinvariant-finite $\Zp\ps{G}$-module. Then we have  an isomorphism
\[ \Ext^1_{\Zp\ps{G}}(M,\Zp\ps{G})\cong \Big( \ilim_U M_U\Big)^\vee.\]
\et

Before giving the proof, we make two remarks.

\br \label{counterex}
Note that the conclusion of theorem is false if we remove the ``systematically coinvariant-finite'' hypothesis. We give a counterexample to illustrate this. Suppose that $d\geq 2$. Let $M=\Zp$ be the module with trivial $G$-action. Then one has $\Ext^1_{\Zp\ps{G}}(M,\Zp\ps{G}) =0$ but $\ilim_U M_U =\Qp$.
\er

\br \label{Gacompare}
It is instructive to specialize our theorem to the context of $G=\Ga\cong\Zp$. In this context, $\Zp\ps{\Ga}$ identifies with the power series ring $\Zp\ps{T}$ in one variable. Denote by $w_n:=w_n(T)$ the polynomial $(1+T)^{p^n}-1$. For a $\Zp\ps{T}$-module, one has a natural identification $M/w_n M \cong M_{\Ga_n}$, where $\Ga_n$ is the subgroup of $\Ga$ of index $p^n$. If $M$ is systematically coinvariant-finite, then the prime ideals dividing $w_n$, $n\geq 0$, are disjoint to the set of prime ideals of height one in the support of $M$. Our Theorem \ref{alg main} therefore coincides with the isomorphism
\[ \Ext^1_{\Zp\ps{\Ga}}(M,\Zp\ps{\Ga}) \cong  \Big( \ilim_n M/w_n\Big)^\vee \]
obtained via the theory of Iwasawa adjoint as in \cite[Proposition 5.5.6]{NSW}.
\er

We split the proof of Theorem \ref{alg main} into a few lemmas. The first of which is the following general observation.

\bl \label{Ext lim}
For every finitely generated $\Zp\ps{G}$-module $M$, there is an isomorphism
\[ \Ext^1_{\Zp\ps{G}}(M,\Zp\ps{G})\cong \plim_U \Ext^1_{\Zp\ps{G}}(M,\Zp[G/U]),\]
where the inverse limit is taken with respect to the canonical projection maps.
\el

\bpf
 The isomorphism
\[ \Hom_{\Zp\ps{G}}(M,\Zp\ps{G})\cong \plim_U \Hom_{\Zp\ps{G}}(M,\Zp[G/U]).\]
plainly holds when $M$ is free of finite rank. But since the ring $\Zp\ps{G}$ is Noetherian, every finitely generated module $M$ has a resolution consisting of finitely generated free $\Zp\ps{G}$-modules. Applying the above isomorphism to this free resolution and taking homology, we obtain the conclusion of the lemma.
\epf

Before continuing, we make another remark.

\br \label{U-Hom}
 Let $U$ be an open (normal) subgroup of $G$ and let $M$ be a $\Zp\ps{G}$-module $M$. A $\Zp\ps{G}$-module homomorphism $M\lra \Zp[G/U]$ naturally factors through $M_U$ to yield a $\Zp\ps{G}$-module homomorphism $M_U\lra \Zp[G/U]$. Furthermore, one sees easily that this latter map can be viewed as a $\Zp[G/U]$-module homomorphism. Conversely, given a $\Zp[G/U]$-module homomorphism $M_U\lra \Zp[G/U]$, by considering the composition $M\lra M_U\lra \Zp[G/U]$, we obtain a $\Zp\ps{G}$-module homomorphism. In conclusion, we have identifications
 \[ \Hom_{\Zp\ps{G}}(M,\Zp[G/U]) = \Hom_{\Zp\ps{G}}(M_U,\Zp[G/U]) =\Hom_{\Zp[G/U]}(M_U,\Zp[G/U]). \]
 These identifications will be frequently utilized in the subsequent discussion of this paper without any further mention.
\er

For the next two lemmas, we shall require the module $M$ to be systematically coinvariant-finite. In particular, the next lemma also makes use of the property of $G$ being commutative.

\bl \label{Ext finite}
Suppose that $M$ is a systematically coinvariant-finite $\Zp\ps{G}$-module. Then we have an isomorphism
\[ \theta_U: \Ext^1_{\Zp\ps{G}}(M,\Zp[G/U])\cong  \Ext^1_{\Zp[G/U]}(M_U,\Zp[G/U]).\]
Furthermore, if $V$ is another open subgroup of $G$ which is contained in $U$, then we have the following commutative diagram
\begin{equation} \label{diag def} \entrymodifiers={!! <0pt, .8ex>+} \SelectTips{eu}{}\xymatrix{
     \Ext^1_{\Zp\ps{G}}(M,\Zp[G/V]) \ar[r]^{\theta_V}_{\sim}\ar[d] &  \Ext^1_{\Zp[G/V]}(M_V,\Zp[G/V]) \ar[d]\\
     \Ext^1_{\Zp\ps{G}}(M,\Zp[G/U]) \ar[r]^{\theta_U}_{\sim} &  \Ext^1_{\Zp[G/U]}(M_U,\Zp[G/U])
     } \end{equation}
     where the vertical maps are induced by the projection $\Zp[G/V]\twoheadrightarrow \Zp[G/U]$.
     \el

\bpf
 Since $M$ is finitely generated over $\Zp\ps{G}$, we can find a short exact sequence
 \[ 0\lra N\lra \Zp\ps{G}^r\lra M\lra 0\]
 of finitely generated $\Zp\ps{G}$-modules. Taking $U$-homology of this short exact sequence, we obtain an exact sequence
\[ 0\lra H_1(U,M)\lra N_U\lra \Zp[G/U]^r\lra M_U \lra 0 \]
of $\Zp[G/U]$-modules. Let $C$ denote the kernel of the map $ \Zp[G/U]^r\lra M_U$. Consider the following diagram
 \begin{equation} \label{hom diagram} \xymatrixcolsep{0.125in}\entrymodifiers={!! <0pt, .85ex>+}\small \SelectTips{eu}{}\xymatrix{
  \Hom_{\Zp[G/U]}(\Zp[G/U]^r,\Zp[G/U]) \ar[r] \ar@{=}[d]&  \Hom_{\Zp[G/U]}(C,\Zp[G/U]) \ar[r] \ar[d]_{\alpha_U}&  \Ext^1_{\Zp[G/U]}(M_U,\Zp[G/U]) \ar[r] \ar@{-->}[d]_{\theta_U}  & 0 \\
\Hom_{\Zp\ps{G}}(\Zp\ps{G}^r,\Zp[G/U]) \ar[r]  &  \Hom_{\Zp\ps{G}}(N ,\Zp[G/U]) \ar[r]  &  \Ext^1_{\Zp\ps{G}}(M,\Zp[G/U]) \ar[r]  & 0 \\
       } \end{equation}
   with exact rows. Here the middle vertical map $\alpha_U$ is given by
\[ \Hom_{\Zp[G/U]}(C,\Zp[G/U])\lra \Hom_{\Zp[G/U]}(N_U,\Zp[G/U]) = \Hom_{\Zp\ps{G}}(N,\Zp[G/U]).\]
It is straightforward to check that the leftmost square is commutative and this in turn induces the rightmost vertical map which is our $\theta_U$. Furthermore, the map $\alpha_U$ is injective with cokernel contained in
\[ \Hom_{\Zp[G/U]}\big(H_1(U,M), \Zp[G/U]\big). \]
Since $U$ is an open subgroup of $G$, it is also isomorphic to $\Zp^d$ as a group. In view that $M_U$ is finite, we also have that $H_1(U,M)$ is finite (cf. \cite[Chap.
IV, Theorem 1]{Se}; one may also consult p.\ 106 in \textit{op. cit.}, where they obtain such a finiteness result for $k\ps{G}$, where $k$ is a field. But the same discussion carries over if $k$ is replaced by $\Zp$.) Since a $\Zp[G/U]$-module homomorphism is also a group homomorphism and $\Zp[G/U]$ has no $p$-torsion, one must have $\Hom_{\Zp[G/U]}\big(H_1(U,M), \Zp[G/U]\big)=0$. Consequently, the map $\alpha_U$ is an isomorphism. From this and the diagram (\ref{hom diagram}), we have that $\theta_U$ is an isomorphism. Finally, one checks easily that the leftmost square in (\ref{hom diagram}) is natural in $U$, which in turn implies that the map $\theta_U$ is natural in $U$. The proof of the lemma is therefore completed. \epf

\bl \label{Ext finite2}
Suppose that $M$ is a systematically coinvariant-finite  $\Zp\ps{G}$-module. Then we have an isomorphism
\[ \psi_U: \Ext^1_{\Zp[G/U]}(M_U,\Zp[G/U])\cong  (M_U)^\vee.\]
Furthermore, if $V$ is another open subgroup of $G$ which is contained in $U$, there is a commutative diagram
\begin{equation} \label{diag Ext2} \entrymodifiers={!! <0pt, .8ex>+} \SelectTips{eu}{}\xymatrix{
     \Ext^1_{\Zp[G/V]}(M_V,\Zp[G/V]) \ar[r]^(.7){\psi_V}_(.7){\sim} \ar[d] &  (M_V)^\vee \ar[d]\\
     \Ext^1_{\Zp[G/U]}(M_U,\Zp[G/U]) \ar[r]^(.7){\psi_U}_(.7){\sim} &  (M_U)^\vee
     } \end{equation}
     where the vertical map on the left is induced by the projection $\Zp[G/V]\twoheadrightarrow \Zp[G/U]$ and the vertical map on the right is induced by the norm map $N_{U/V}: M_U\lra M_V$.
     \el

\bpf
Since $M_U$ is finite, it is annihilated by $p^t$ for some large enough $t$. On the other hand, multiplication by $p^t$ induces an automorphism on $\Zp[G/U]\ot_{\Zp}\Qp$. Hence it follows that \[\Ext^i_{\Zp[G/U]}(M_U,\Zp[G/U]\ot_{\Zp}\Qp)=0\] for every $i\geq 0$. Taking this into account, upon applying $\Hom_{\Zp[G/U]}(M_U,-)$ to the short exact sequence
\[  0\lra \Zp[G/U] \lra \Zp[G/U]\ot_{\Zp}\Qp \lra \Zp[G/U]\ot_{\Zp}\Qp/\Zp\lra 0, \]
we see that the connecting morphism
\[  \partial_U: \Hom_{\Zp[G/U]}(M_U, \Zp[G/U]\ot_{\Zp}\Qp/\Zp)\lra \Ext^1_{\Zp[G/U]}(M_U,\Zp[G/U]) \]
is an isomorphism which is natural in $U$. It therefore remains to show the existence of an isomorphism
\[ \eta_U: \Hom_{\Zp[G/U]}(M_U, \Zp[G/U]\ot_{\Zp}\Qp/\Zp)\cong (M_U)^\vee \]
for every open subgroup $U$ of $G$ with the property that if $V\subseteq U$, there is a commutative diagram
\begin{equation} \label{diag Hom vee} \entrymodifiers={!! <0pt, .8ex>+} \SelectTips{eu}{}\xymatrix{
     \Hom_{\Zp[G/V]}(M_V,\Zp[G/V]\ot_{\Zp}\Qp/\Zp) \ar[r]^(.77){\eta_V}_(.77){\sim} \ar[d] &  (M_V)^\vee \ar[d]\\
     \Hom_{\Zp[G/U]}(M_U,\Zp[G/U]\ot_{\Zp}\Qp/\Zp) \ar[r]^(.77){\eta_U}_(.77){\sim} &  (M_U)^\vee
     } \end{equation}
     where the left vertical map is induced by the projection $\Zp[G/V]\twoheadrightarrow \Zp[G/U]$ and the right vertical map is induced by the norm map $N_{U/V}: M_U\lra M_V$.
     To simplify notation, we shall write $\Qp/\Zp[G/U] = \Zp[G/U]\ot_{\Zp}\Qp/\Zp$.
Now, for each $f\in \Hom_{\Zp[G/U]}(M_U, \Qp/\Zp[G/U])$ and $x\in M_U$, we write
\[ f(x) =\sum_{\sigma\in G/U}f_{\sigma}(x)\sigma,\]
where $f_{\sigma}(x)\in\Qp/\Zp$. We then define the map
\[ \eta_U: \Hom_{\Zp[G/U]}(M_U, \Qp/\Zp[G/U])\lra (M_U)^\vee \]
by sending $f\mapsto f_1$, where ``$1$" denotes the identity of the group $G/U$. Let $\tau\in G/U$. Since $f$ is a $\Zp[G/U]$-homomorphism, we have the identity $\tau f(x) = f(\tau x)$ which in turn yields
\[ f_\sigma(\tau x) = f_{\tau^{-1}\sigma}(x) \]
for every $\sigma\in G/U$. In particular, we have
\[ f_1(\tau x) = f_{\tau^{-1}}(x) \]
for every $\tau\in G/U$. Hence $f$ is uniquely determined by $f_1$, and therefore, $\eta_U$ is an isomorphism. It remains to show that $\eta_U$ is natural in $U$. Let $V$ be an open subgroup of $G$ which is contained in $U$. The projection $\Zp[G/V]\lra \Zp[G/U]$ induces a map
\[ \rho_{VU}:\Hom_{\Zp[G/V]}(M_V, \Qp/\Zp[G/V])\lra  \Hom_{\Zp[G/V]}(M_V, \Qp/\Zp[G/U]) = \Hom_{\Zp[G/U]}(M_U, \Qp/\Zp[G/U]).\]
Let $x\in M_U$ and $h\in \Hom_{\Zp[G/V]}(M_V, \Zp[G/V]\ot_{\Zp}\Qp/\Zp)$. By abuse of notation, we write $x\in M_V$ for a preimage of $x$ under natural surjection $M_V\lra (M_V)_{U/V}=M_U$. Let $\sigma_1=1, \sigma_2,...,\sigma_r \in G/V$ be a complete set of coset representatives of $U/V$ in $G/V$. Then one has
\[ h(x) = \sum_{\tau\in U/V}\sum_{i=1}^r h_{\tau\sigma_i}(x)\tau\sigma_i.\]
From which, we have
\[ \rho_{VU}(h) (x) = \rho_{VU}(h(x)) = \sum_{i=1}^r\Big(\sum_{\tau\in U/V}h_{\tau\sigma_i}(x)\Big)\sigma_i,  \]
which in turn implies that
\[ \rho_{VU}(h)_1 (x) =\sum_{\tau\in U/V}h_{\tau}(x) = \sum_{\tau\in U/V} h_1(\tau^{-1}x) = h_1\left(\sum_{\tau\in U/V} \tau^{-1} x \right) =h_1\big(N_{U/V}(x)\big). \]
This establishes the commutativity of the diagram (\ref{diag Hom vee}). The proof is therefore completed.
\epf

Theorem \ref{alg main} is now a consequence of a combination of Lemmas \ref{Ext lim}, \ref{Ext finite} and \ref{Ext finite2}. We end this section with two useful corollaries. Recall that for a $\Zp\ps{G}$-module $M$, $M^\iota$ denotes the $\Zp\ps{G}$-module which is the same underlying $\Zp$-module $M$ with $G$-action given by
\[ g \cdot_{\iota} x = g^{-1}x, \quad g \in G,~ x\in M.\]

\bc \label{dual pseudo}
Let $M$ be a systematically coinvariant-finite $\Zp\ps{G}$-module. Then one has a pseudo-isomorphism
\[ M^\iota\cong \Big( \ilim_U M_U\Big)^\vee.\]
of $\Zp\ps{G}$-modules. \ec

\bpf
 Since $\Ext^1_{\Zp\ps{G}}(M, \Zp\ps{G})\sim M^\iota$ (cf. \cite[Proposition 8]{PR}), the conclusion of the corollary follows from this and Theorem \ref{alg main}.
\epf

\bc \label{pseudo and zero}
Let $M$ be a systematically coinvariant-finite $\Zp\ps{G}$-module. Then $M$ is pseudo-null over $\Zp\ps{G}$ if and only if $\ilim_U M_U =0$. \ec

\bpf
By virtue of Lemma \ref{finite fg tor}, we already know that $M$ is finitely generated torsion over $\Zp\ps{G}$. Therefore, for $M$ to be pseudo-null over $\Zp\ps{G}$, it is equivalent to having $\Ext^1_{\Zp\ps{G}}(M,\Zp\ps{G})=0$. In view of Theorem \ref{alg main}, this is the same as saying that $\ilim_U M_U =0$.
\epf

\section{Arithmetic} \label{Arithmetic preliminaries}

Let $F$ be a number field. In the event $p=2$, we assume further that $F$ has no real places. For a ring $R$ with identity, write $K_n(R)$ for the algebraic $K$-groups of $R$ in the sense of Quillen \cite{Qui73a} (also see \cite{Kol, WeiKbook}). Let $i\geq 2$. We let $\Op_F$ denote the ring of integers of $F$. By the fundamental results of Quillen \cite{Qui73b} and Borel \cite{Bo}, the group $K_{2i-2}(\Op_F)$ is finite for each $i\geq 2$. For a finite extension $L$ of $F$, we have a map
 \[ \jmath_{F/L}:  K_{2i-2}(\Op_F) \lra K_{2i-2}(\Op_L)\]
induced by the inclusion $\Op_F\lra \Op_L$ via functionality. On the other hand, we have the transfer map
\[ \Tr_{L/F}:  K_{2i-2}(\Op_L) \lra K_{2i-2}(\Op_F).\]

Let $F_\infty$ be a $\Zp^d$-extension of $F$, where $d\geq 1$. The Galois group $\Gal(F_\infty/F)$ will always be denoted by $G$. We then consider the following direct limit and inverse limit
\[ \ilim_L K_{2i-2}(\Op_L)[p^\infty] \quad\mbox{and}\quad \plim_L K_{2i-2}(\Op_L)[p^\infty],\]
where the transition maps for the direct limit (resp., the inverse limit) are given by $\jmath_{L/L'}$ (resp., $\Tr_{L/L'}$). For the direct limit, we shall sometimes write
\[ K_{2i-2}(\Op_{F_\infty})_p: = \ilim_L K_{2i-2}(\Op_L)[p^\infty].\]

The following is the main theorem of this paper.

\bt \label{main theorem}
For $i\geq 2$, there is a pseudo-isomorphism
\[ \Big(\plim_L K_{2i-2}(\Op_L)[p^\infty]\Big)^\iota \sim \Big( K_{2i-2}(\Op_{F_\infty})_p\Big)^\vee.\]
of $\Zp\ps{G}$-modules. \et

As seen from the discussion in the introduction, it suffices to show the following assertion.

\bp \label{main prop}
The $\Zp\ps{G}$-module $\plim_L K_{2i-2}(\Op_L)[p^\infty]$ is systematically coinvariant-finite such that for every open subgroup $U$ of $G$, there is an isomorphism
      \[t_U: \left(\plim_L K_{2i-2}(\Op_L)[p^\infty]\right)_U \cong K_{2i-2}(\Op_{L_U})[p^\infty],\] where $L_U$ is the fixed field of $U$. Furthermore, if $V$ is an open subgroup of $G$ contained in $U$ with fixed field $L_V$, we then have the following commutative diagram \begin{equation} \label{K commute} \entrymodifiers={!! <0pt, .8ex>+} \SelectTips{eu}{}\xymatrix{
     \left(\plim_L K_{2i-2}(\Op_L)[p^\infty]\right)_U \ar[r]_(0.55){\sim}\ar[d]^{N_{U/V}} &  K_{2i-2}(\Op_{L_V})[p^\infty] \ar[d]^{\jmath_{L_U/L_V}}\\
     \left(\plim_L K_{2i-2}(\Op_L)[p^\infty]\right)_V \ar[r]_(0.55){\sim} &  K_{2i-2}(\Op_{L_U})[p^\infty]}
     \end{equation}
     \ep

Theorem \ref{main theorem} will then follow immediately from a combination of Corollary \ref{dual pseudo} and Proposition \ref{main prop}. Furthermore, combining Corollary \ref{pseudo and zero} with Proposition \ref{main prop} yields the following corollary.

\bc \label{pseudo K=0}
Retain the above notation. Then $\plim_L K_{2i-2}(\Op_L)[p^\infty]$ is pseudo-null over $\Zp\ps{G}$ if and only $\ilim_L K_{2i-2}(\Op_L)[p^\infty]=0$.
\ec

The remainder of the section will be devoted to the verification of Proposition \ref{main prop}. Throughout, we shall let $S$ denote the set of primes of $F$ consisting of those above $p$ and the infinite primes. Write $F_S$ for the maximal algebraic extension of $F$ unramified outside $S$. Denoting by $\mu_{p^n}$ the cyclic group generated by a primitive $p^n$-root of unity, we then write $\mu_{p^\infty}$ for the direct limit of the groups $\mu_{p^n}$. These have natural $G_S(F)$-module structures. The action of $G_S(F)$ on $\mu_{p^\infty}$ induces a continuous character
\[\chi: G_S(F) \lra \mathrm{Aut}(\mu_{p^\infty}) \cong \Zp^{\times}.\]
For a discrete or compact $G_S(F)$-module $X$, we shall write $X(i)$ for the $G_S(F)$-module which is $X$ as a $\Zp$-module but with a $G_S(F)$-action given by
\[ \sigma\cdot x = \chi(\sigma)^i\sigma x,  \]
where the action on the right is the original action of $G_S(F)$ on $X$.
Plainly, we have $X(0)=X$ and $\mu_{p^{\infty}} \cong \Qp/\Zp(1)$. One can also check directly that
\[X(i+j) \cong \big(X(i)\big)(j).\]

The key approach towards proving Proposition \ref{main prop} is via cohomology. For this, we need to recall the works of \cite{Sou}, Rost and Voevodsky \cite{Vo} which allows us to translate our problem into a cohomological one. In \cite{Sou}, Soul\'e connected the $K$-groups with continuous cohomology groups via the $p$-adic Chern class maps
\[ \mathrm{ch}_{i}^{F}: K_{2i-2}(\Op_F)[p^\infty] \cong K_{2i-2}(\Op_F)\ot \Zp \lra H^2\left(G_S(F), \Zp(i)\right)\]
for $i\geq 2$. For the precise definition of these maps, we refer readers to loc.\ cit.
Soul\'e has proved that these maps are surjective (see \cite[Th\'eor\`{e}me 6(iii)]{Sou}). Thanks to the deep work of Rost and Voevodsky \cite{Vo} (also see \cite{Wei09}), we now know that these maps are isomorphisms.

Now, if $L$ is a finite extension of $F$ contained in $F_S$, we shall write $G_S(L)$ for the Galois group $\Gal(F_S/L)$. Then one has the following commutative diagrams (see \cite[Chap. III]{Sou})
\begin{equation} \label{tr res}\entrymodifiers={!! <0pt, .8ex>+} \SelectTips{eu}{}
\xymatrix{
     K_{2i-2}(\Op_F)[p^\infty]\ar[r]^{\mathrm{ch}_i^F} \ar[d]^{\jmath_{L/F}} &  H^2\left(G_S(F), \Zp(i)\right) \ar[d]^{\mathrm{res}}\\
     K_{2i-2}(\Op_L)[p^\infty] \ar[r]^{\mathrm{ch}_i^L} &  H^2\left(G_S(L), \Zp(i)\right)} \end{equation}
\begin{equation} \label{N cor}\entrymodifiers={!! <0pt, .8ex>+} \SelectTips{eu}{}
\xymatrix{
      K_{2i-2}(\Op_L)[p^\infty] \ar[r]^{\mathrm{ch}_i^L} \ar[d]^{\Tr_{L/F}} &  H^2\left(G_S(L), \Zp(i)\right) \ar[d]^{\mathrm{cor}} \\
      K_{2i-2}(\Op_F)[p^\infty] \ar[r]^{\mathrm{ch}_i^F} &  H^2\left(G_S(F), \Zp(i)\right)}
      \end{equation}

We turn back to our $\Zp^d$-extension $F_{\infty}$. Note that by \cite[Theorem 1]{Iw73}, the extension $F_\infty$ is contained in $F_S$. Hence it makes sense to speak of $G_S(L)=\Gal(F_S/L)$. We then define the Iwasawa cohomology group $H^k_{\Iw}(F_\infty/F, \Zp(i))$ to be
\[H^k_{\Iw}(F_\infty/F, \Zp(i)):= H^k_{\Iw, S}(F_\infty/F, \Zp(i))=\plim_L H^k(G_S(L),\Zp(i)),\]
where the inverse limit is taken over all the finite extensions $L$ of $F$ contained
in $F_\infty$ and with respect to the corestriction maps. It can be shown that these cohomology groups are finitely generated over $\Zp\ps{G}$ (for instance, see \cite[Proposition 4.1.3]{LimSh}). From the commutative diagram (\ref{N cor}), we obtain the following relation between the inverse limit of $K$-groups and the second Iwasawa cohomology groups.

\bl \label{K2 Iw coh}
Suppose that $i\geq 2$. Then there is an isomorphism
\[ \plim_L K_{2i-2}(\Op_L)[p^\infty] \cong  H^2_{\Iw, S}(F_\infty/F,\Zp(i))\] of $\Zp\ps{G}$-modules.
\el

We now recall the following version of Tate's descent spectral sequence for Iwasawa cohomology groups. This will allow us to relate the coinvariant of the Iwasawa cohomology group with the intermediate cohomology groups.

\bp \label{gal descent}
Let $U$ be an open normal subgroup of $G=\Gal(F_\infty/F)$ and write $L$ for the fixed field of $U$. Then we have a homological spectral sequence
 $$ H_r\big(U, H^{-s}_{\Iw}(F_\infty/F, \Zp(i))\big)\Longrightarrow H^{-r-s}\big(G_S(L), \Zp(i)\big). $$
 In particular, we have an isomorphism
\[ H^2_{\Iw}\big(F_\infty/F, \Zp(i)\big)_U \cong H^2\big(G_S(L), \Zp(i)\big).\]
induced by the corestriction map on cohomology.
\ep

\bpf
Had $F_\infty/F$ being a finite extension, this is essentially the Tate spectral sequence (for instance, see \cite[Theorem 2.5.3]{NSW}). In the general context of the proposition, this follows from the work of Nekov\v{a}\'r \cite[Proposition 8.4.8.1]{Ne}.
 The final isomorphism in the proposition follows from reading off the initial $(0,-2)$-term of the spectral sequence.
\epf

In view of the preceding proposition and commutative diagram (\ref{tr res}), for the verification of the commutativity of (\ref{K commute}), one is reduced to proving the following lemma.

\bl \label{N cor Iw}
For open subgroups $U,V$ of $G$ such that $V\subseteq U$, let $L$ (resp., $K$) be the fixed field of $U$ (resp., $V$). Then we have the following commutative diagram
\[ \entrymodifiers={!! <0pt, .8ex>+} \SelectTips{eu}{}\xymatrix{
      H^2_{\Iw}\big(F_\infty/F, \Zp(i)\big)_U \ar[r]^{\mathrm{cor}}_{\sim}\ar[d]_{N_{U/V}} &  H^2\big(G_S(L), \Zp(i)\big) \ar[d]_{\mathrm{res}}\\
      H^2_{\Iw}\big(F_\infty/F, \Zp(i)\big)_V  \ar[r]^{\mathrm{cor}}_{\sim} &  H^2\big(G_S(K), \Zp(i)\big)
     }\]
\el

\bpf
Let $W$ be any open subgroup of $G$ contained in $V$. Write $E=E_W$ for the fixed field of $W$. By either appealing to the double coset formula (cf. \cite[Proposition 1.5.11]{NSW}) or a direct verification, one has a commutative diagram
\[ \entrymodifiers={!! <0pt, .8ex>+} \SelectTips{eu}{}\xymatrix{
     H^2\big(G_S(E), \Zp(i)\big)\ar[r]^{\mathrm{cor}}\ar[d]_{N_{U/V}} &  H^2\big(G_S(L), \Zp(i)\big) \ar[d]_{\mathrm{res}}\\
     H^2\big(G_S(E), \Zp(i)\big) \ar[r]^{\mathrm{cor}} &  H^2\big(G_S(K), \Zp(i)\big)
     }\]
Varying $W$, we obtain a commutative diagram
\[ \entrymodifiers={!! <0pt, .8ex>+} \SelectTips{eu}{}\xymatrix{
    H^2_{\Iw}\big(F_\infty/F, \Zp(i)\big)\ar[r]^{\mathrm{cor}}\ar[d]_{N_{U/V}} &  H^2\big(G_S(L), \Zp(i)\big) \ar[d]_{\mathrm{res}}\\
     H^2_{\Iw}\big(F_\infty/F, \Zp(i)\big) \ar[r]^{\mathrm{cor}} &  H^2\big(G_S(K), \Zp(i)\big)
     }\]
     By Proposition \ref{gal descent}, the horizontal maps factor through $H^2_{\Iw}\big(F_\infty/F, \Zp(i)\big)_U$ and $H^2_{\Iw}\big(F_\infty/F, \Zp(i)\big)_V$ respectively. From which, we obtain the required diagram of the lemma.
\epf

We may now conclude the section.

\bpf[Proof of Proposition \ref{main prop}]
The isomorphism $t_U$ is a consequence of Lemma \ref{K2 Iw coh} and Proposition \ref{gal descent}. On the other hand, the commutativity of diagram (\ref{K commute}) follows from a combination of Proposition \ref{gal descent}, Lemma \ref{N cor Iw} and diagram (\ref{tr res}).
\epf

\section{Some further remarks} \label{examples and remark}

A conjecture of
Schneider (cf.\ \cite[p.\ 192]{Sch79}) asserted that the cohomology group $H^2(G_S(F),\Zp(i))$ is finite for $i\leq 0$. Granted this conjecture, the argument in the preceding section carries over to yield a similar sort of result for these cohomology groups. We will just record two situations, where the conjecture of Schneider is known to hold.

\bt
Suppose that we are in either of the following situations.
\begin{enumerate}
  \item[(a)] The number field $F$ is totally real (and so $p\geq 3$ by our standing assumption), $F_\infty$ is the cyclotomic $\Zp$-extension of $F$ and $i$ is a negative odd integer.
  \item[(b)] The number field $F$ is an imaginary quadratic field, $F_\infty$ is a $\Zp^d$-extension of $F$ (note that $d=1$ or 2) and $i=0$.
\end{enumerate}
Then we have a pseudo-isomorphism
\[ \Big(\plim_L H^2(G_S(L),\Zp(i))\Big)^{\iota}\sim \Big(\ilim_L H^2(G_S(L),\Zp(i))\Big)^{\vee}\]
of $\Zp\ps{\Gal(F_\infty/F)}$-modules.
\et

\bpf
The proof is similar as before. We simply mention that in case (a), the finiteness of $H^2(G_S(L),\Zp(i))$ is a consequence of \cite[Proposition 3.8]{Nic}. For case (b), since each intermediate extension $L$ of $F$ contained in $F_\infty$ is abelian over $F$, Brumer's theorem (cf.\ \cite[Theorem 10.3.16]{NSW}) applies telling us that
$H^2(G_S(L),\Zp)$ is finite.
\epf

\section{Some classes of examples} \label{Examples}

We end giving some classes of examples, where $\ilim_L K_{2i-2}(\Op_L)[p^\infty]$ can be either zero or not. We first show that one can construct many examples of nonzero $\ilim_L K_{2i-2}(\Op_L)[p^\infty]$. Recall that for every $\Zp\ps{G}$-module $M$, there is a $\Zp\ps{G}$-homomorphism
\[ M[p^\infty] \lra \bigoplus_{i=1}^s\Zp\ps{G}/p^{\alpha_i} \]
with pseudo-null kernel and cokernel. The $\mu_G$-invariant of $M$ is then defined to be $\sum_{i}^{s}\alpha_i$. Note that if $M$ (and hence $M[p^\infty]$) is pseudo-null over $\Zp\ps{G}$, then its $\mu_G$-invariant is trivial.

\bp
Let $i\geq 2$ and let $G = \Zp^d$, where $d\geq 1$. Suppose that $p$ is a prime such that
$p > 2d+1$. Then there exist infinitely many pairs $(F, F_\infty)$, where $F$ is a finite cyclic extension of $\Q(\mu_p)$ and
$F_\infty$ is a $\Zp^d$-extension of $F$ such that
$K_{2i-2}(\Op_{F_\infty})_p\neq 0$.
\ep

\bpf
Indeed, under the hypothesis of the result, it has been shown that there exist infinitely many pairs $(F, F_\infty)$, where $F$ is a finite cyclic extension of $\Q(\mu_p)$ and
$F_\infty$ is a $\Zp^d$-extension of $F$ such that $H^2_{\Iw}(F_\infty/F,\Zp(i))$ has nontrivial $\mu_G$-invariant (cf.\ \cite[Proposition 5.2.2]{LimKgroups}). By the remark before the proposition, the module $H^2_{\Iw}(F_\infty/F,\Zp(i))$ is therefore not pseudo-null over $\Zp\ps{G}$. Consequently, it follows from Lemma \ref{K2 Iw coh}
and Corollary \ref{pseudo K=0} that $K_{2i-2}^p(\Op_{F_\infty})\neq 0$. \epf

\bp \label{Sharifi}
 Let $F=\Q(\mu_p)$, where $p$ is an irregular prime $<1000$. If $F_\infty$ is the compositum of all $\Zp$-extensions of $F$, then $K_{2i-2}(\Op_{F_\infty})_p= 0$ for every $i\geq 2$.
\ep

\bpf
Let $L_\infty$ be the maximal unramified abelian pro-$p$ extension of $F_{\infty}$, and let $L_\infty'$ be the maximal subextension of $L_\infty$ in which every
prime of $F_{\infty}$ above $p$ splits completely. Sharifi has shown that the Greenberg conjecture is valid under the hypothesis of the proposition (cf. \cite[Theorem 1.3]{Sh08}). In other words, $\Gal(L_\infty/F_\infty)$ is pseudo-null over $\Zp\ps{G}$. Since $\Gal(L'_\infty/F_\infty)$ is a quotient of  $\Gal(L_\infty/F_\infty)$, it is also pseudo-null over $\Zp\ps{G}$. Now, by the Poitou-Tate sequence, we have an exact sequence
\[ 0\lra \Gal(L'_\infty/F_\infty)\lra H_{\Iw}^2(F_\infty/F,\Zp(1))\lra \plim_{L} \bigoplus_{w_L\in S(L)} H^2(L_{w_L},\Zp(1)), \]
where $S(L)$ is the set of primes of $L$ above $p$.
Since the decomposition group of $F_\infty/F$ at the prime of $F$ above $p$ has dimension $\geq 2$ (cf. \cite[Th\'{e}or\`{e}me 3.2]{LaNQD}), we may apply a similar argument to that in \cite[Lemma 5.3]{LimFine} to conclude that $ \plim_{L} \bigoplus_{w_L\in S(L)} H^2(L_{w_L},\Zp(1))$ is pseudo-null over $\Zp\ps{G}$, and whence, $H_{\Iw}^2(F_\infty/F,\Zp(1))$ is pseudo-null over $\Zp\ps{G}$. On the other hand, as $F_\infty$ contains $\mu_{p^{\infty}}$, an application of \cite[Lemma 2.5.1(c)]{Sh22} tells us that
\[ H_{\Iw}^2(F_\infty/F,\Zp(i))\cong H_{\Iw}^2(F_\infty/F,\Zp(1))\ot\Zp(i-1) \]
for every $i\geq 2$. Thus, it follows that $H_{\Iw}^2(F_\infty/F,\Zp(i))$ is pseudo-null over $\Zp\ps{G}$ for every $i\geq 2$. Combining this latter observation with Lemma \ref{K2 Iw coh} and Corollary \ref{pseudo K=0}, we obtain the conclusion of the proposition. \epf


\footnotesize
\begin{thebibliography}{00}
\bibitem{BH} P. N. Balister; S. Howson, Notes on Nakayama's
lemma for compact $\La$-modules, Asian J. Math. 1(2) (1997) 224--229.

\bibitem{Ba} A. Bandini, Greenberg's conjecture and capitulation
in $\Zp^d$-extensions, J. Number Theory 122 (2007) 121--134.

\bibitem{Bo} A. Borel, Stable real cohomology of arithmetic groups, Ann. Sci. \'Ecole Norm. Sup. (4) 7 (1974), 235--272.


\bibitem{G01} R. Greenberg,  Iwasawa theory--past and present, in: Class field theory--its centenary and prospect (Tokyo, 1998), 335--385, Adv. Stud. Pure Math., 30, Math. Soc. Japan, Tokyo, 2001.

\bibitem{Iw73} K. Iwasawa, On ${\bf Z}\sb{l}$-extensions of algebraic number fields, Ann. of Math. (2) 98 (1973), 246--326.




\bibitem{Kol} M. Kolster, $K$-theory and arithmetic. Contemporary developments in algebraic $K$-theory, 191--258, ICTP Lect. Notes, XV, Abdus Salam Int. Cent. Theoret. Phys., Trieste, 2004.

\bibitem{LT} K. F. Lai; K.-S. Tan,  A generalized Iwasawa's theorem and its application, Res. Math. Sci. 8 (2021), no. 2, Paper No. 20, 18 pp.

\bibitem{LaNQD} A. Lannuzel; T. Nguyen Quang Do, Conjectures de Greenberg et extensions pro-$p$-libres d'un corps de nombres, Manuscripta Math. 102 (2000), no. 2, 187--209.

\bibitem{LimFine}  M. F. Lim, Notes on the fine Selmer groups, Asian J. Math. 21(2) (2017) 337--362.

\bibitem{LimKgroups} M. F. Lim,  On the growth of even $K$-groups of rings of integers in $p$-adic Lie extensions, accepted for publication in Israel J. Math., arXiv:2009.01477[math.NT].

\bibitem{LimSh} M. F. Lim; R. Sharifi,
Nekov\'a\v{r} duality over $p$-adic Lie extensions of global fields, Doc. Math. 18 (2013), 621--678.


\bibitem{Ne} J. Nekov\'a\v{r}, Selmer complexes, Ast\'erisque No. 310 (2006), viii+559 pp.

\bibitem{NSW} J. Neukirch; A. Schmidt; K. Wingberg,
Cohomology of Number Fields. 2nd edn., Grundlehren Math.
Wiss. 323 (Springer-Verlag, Berlin, 2008).

\bibitem{Nic} A. Nickel, Annihilating wild kernels, Doc. Math. 24 (2019), 2381--2422.

\bibitem{PR} B. Perrin-Riou, Arithm\'etique des courbes elliptiques et théorie d'Iwasawa, M\'em. Soc. Math. France (N.S.) No. 17 (1984), 130 pp.

\bibitem{Sch79} P. Schneider, \"{U}ber gewisse Galoiscohomologiegruppen, Math. Z. 168 (1979), no. 2, 181--205.

\bibitem{Se} J. -P. Serre, Local Algebra. Springer Monographs in Mathematics, Spring-Verlag, Berlin, 2000

\bibitem{Sh08} R. Sharifi, On Galois groups of unramified pro-$p$ extensions, Math. Ann. 342 (2008), no. 2, 297--308.

\bibitem{Sh22}    R. Sharifi, Reciprocity maps with restricted ramification, arXiv:1609.0361v2[math.NT].

\bibitem{Sou} C. Soul\'e, $K$-th\'eorie des anneaux d'entiers de corps de nombres et cohomologie \'etale, Invent. Math. 55 (1979), no. 3, 251--295.

\bibitem{Qui73a} D. Quillen, Higher algebraic $K$-theory. I. Algebraic $K$-theory, I: Higher $K$-theories (Proc. Conf., Battelle Memorial Inst., Seattle, Wash., 1972), pp. 85--147. Lecture Notes in Math., Vol. 341, Springer, Berlin 1973.

\bibitem{Qui73b} D. Quillen, Finite generation of the groups $K\sb{i}$ of rings of algebraic integers. Algebraic $K$-theory, I: Higher $K$-theories (Proc. Conf., Battelle Memorial Inst., Seattle, Wash., 1972), pp. 179--198. Lecture Notes in Math., Vol. 341, Springer, Berlin, 1973.

\bibitem{Vau} D. Vauclair, Sur la dualit\'e et la descente d'Iwasawa, Ann. Inst. Fourier Grenob. 59(2) (2009) 691--767.


\bibitem{Vo} V. Voevodsky, On motivic cohomology with $\mathbf Z/l$-coefficients, Ann. of Math. (2) 174 (2011), no. 1, 401--438.

\bibitem{Wei09} C. Weibel, The norm residue isomorphism theorem, J. Topol. 2 (2009), no. 2, 346--372.

\bibitem{WeiKbook} C. Weibel, The $K$-book. An introduction to algebraic $K$-theory. Graduate Studies in Mathematics, 145. American Mathematical Society, Providence, RI, 2013. xii+618 pp.

\end{thebibliography}
\end{document}